\documentclass[10pt]{article}
\usepackage{amssymb}

\newtheorem{lemma}{Lemma}[section]

\newtheorem{theorem}{Theorem}[section]

\begin{document}

\title{Topological invariants of three-manifolds from $U_{q}(osp(1|2n))$}
\author{ {Sacha C. Blumen\footnote{e-mail: sachab@maths.usyd.edu.au}} \\
{\small {\em School of Mathematics and Statistics,}} \\ 
{\small {\em University of Sydney, NSW 2006, Australia}} }
\maketitle

\begin{abstract} 
We create Reshetikhin-Turaev topological invariants of closed orientable three-manifolds from the
quantum supergroup $U_{q}(osp(1|2n))$ at certain even roots of unity.  To construct the invariants we
develop tensor product theorems for finite dimensional modules of $U_{q}(osp(1|2n))$ at roots of unity.
\end{abstract}

\begin{section}{Introduction}

Topological invariants of closed orientable three-manifolds
may be constructed from modular or quasimodular Hopf algebras \cite{rt,tw}.
Reshetikhin and Turaev's construction using modular Hopf algebras relies upon
several theorems relating framed links in $S^{3}$ to closed orientable
three-manifolds.  The Lickorish-Wallace theorem states that each
framed link in $S^{3}$ determines a closed, orientable 3-manifold and that every
such 3-manifold is obtainable by performing surgery upon a framed link in
$S^{3}$.  Kirby-Craggs, and Fenn and Rourke showed that homeomorphism
classes of closed orientable three-manifolds may be generated by performing surgery
upon elements of equivalence classes of framed links in $S^{3}$, 
where the equivalence relations are generated by the Kirby moves.
By taking such combinations of isotopy invariants of links in $S^{3}$
as to render them unchanged under the Kirby moves one obtains a
topological invariant of 3-manifolds.

RT took invariants of isotopy derived from the quantum group $U_{q}(sl_{2})$ at
even roots of unity.  Their method was adapted for the
quantum algebras related to the $A_{n}, B_{n}, C_{n}, D_{n}$ Lie algebras at all
roots of unity \cite{tw,z1}, and the exceptional quantum algebras
and quantum superalgebras $U_{q}(osp(1|2))$ and $U_{q}(gl(2|1))$ at odd roots of
unity \cite{z2}.  Here we create invariants from $U_{q}(osp(1|2n))$ at
$q=\exp{(2\pi i/N)}$ with $N=2(2k+1)$.

Resethikhin-Turaev invariants can be constructed from a class of Hopf (super)algebras more general than
(quasi)modular Hopf (super)algebras.
For a Hopf (super)algebra $A$, invariants may be constructed if $A$ has the following
properties 
(where we take the quantum superdimension and quantum supertrace if $A$ is a quantum superalgebra)
\begin{itemize}
\item[i)] there exists a finite collection of mutually non-isomorphic left
$A$-modules $\{ V_{\lambda} \}$ where $\lambda$ ranges over some index set $I$
such that dim$(V_{\lambda}) < \infty$ and $dim_{q}(V_{\lambda}) \neq 0$, $\forall \lambda \in I$,
\item[ii)] for any finite collection of $A$-modules $V_{\lambda_{1}}, V_{\lambda_{2}},
\ldots, V_{\lambda_{s}}$ such that $\lambda_{i} \in I$ for all $i$,
$$V_{\lambda_{1}} \otimes V_{\lambda_{2}} \otimes \cdots \otimes
V_{\lambda_{s}} = {\cal V} \oplus {\cal Z}$$ where ${\cal V} =
\bigoplus_{\lambda \in I}(V_{\lambda})^{\oplus m(\lambda)}$, $m(\lambda)
\geq 0$ is the multiplicity of $V_{\lambda}$ in the direct sum, 
and ${\cal Z}$ is a possibly empty $A$-module with zero quantum dimension,
\item[iii)] for each $V_{\lambda}$ there is a dual module
$(V_{\lambda})^{\dagger} \cong (V_{\lambda^{*}})$ such that $\lambda^{*} \in I$
and there exists a distinguished module $V_{0}$ such that $(V_{0})^{\dagger}
\cong V_{0}$,
\item[iv)] the central element $\delta = v - \sum_{\lambda \in I}d_{\lambda}
\chi_{\lambda}(v^{-1})C_{\lambda}$ vanishes upon acting on any $V_{\mu}$ where $\mu \in \Lambda_{N}^{+}$.  
Here $\chi_{\zeta}(v) = q^{-(\zeta,\zeta+2\rho)}$ and $C_{\lambda} =
tr_{\lambda}[({\textrm {id}} \otimes \pi)(1 \otimes q^{2h_{\rho}})R^{T}R]$. $R$ is the universal $R$-matrix and $R^{T}
= P.R.P$ where $P$ is the permutation operator.  
$\{d_{\lambda} \}$ is a collection of complex valued constants such that at least
one $d_{\lambda}$ is non-zero.  This condition ensures that combinations of
isotopy invariants of links are unchanged under some of the Kirby moves, and
\item[v)] the sum $z = \sum_{\lambda \in I}d_{\lambda} q^{-(\lambda,\lambda +
2\rho)}dim_{q}(V_{\lambda})$ is non-zero.
\end{itemize}

\end{section}

\begin{section}{$U_{q}(osp(1|2n))$ at roots of unity and its finite dimensional modules}
The quantum superalgebra $U_{q}(osp(1|2n))$ at roots of unity is not quasi-triangular.  However, the quantum
superalgebra now has a class of central elements which do not exist at generic $q$.  These central elements
generate an ideal of $U_{q}(osp(1|2n))$, which is also a two-sided co-ideal.

The quotient of $U_{q}(osp(1|2n))$ by this ideal turns out to be a quasi-triangular Hopf superalgebra.
Let us denote this algebra by $U_{q}^{(N)}(osp(1|2n))$, the details of which may be found in \cite{z5}. 

Associated with $osp(1|2n)$ there is a euclidean space $H^{*}$ which has a basis of vectors
$\{ \epsilon_{i} | 1 \leq i \leq n \}$ and a bilinear form such that $(\epsilon_{i},\epsilon_{j}) =
\delta_{i,j}$.  The even (resp. odd) positive roots are 
$\Phi_{0}^{+} = \{\epsilon_{i} \pm \epsilon_{j}, 2\epsilon_{k} | 1 \leq i < j \leq n, 1 \leq k \leq n \}$ 
(resp.  $\Phi_{1}^{+} = \{ \epsilon_{i} | 1 \leq i \leq n \}$).  
Define $2\rho = \sum_{\alpha \in \Phi_{0}^{+}}\alpha - \sum_{\beta \in \Phi_{1}^{+}}\beta$, 
$\mathbb{Z} = \{ \ldots, -1, 0, 1, \ldots \}$, 
$\mathbb{Z}_{+} = \{0,1,2,\ldots \}$ and $\mathbb{Z}_{n} = \{0,1,2,\ldots,n-1\}$.  
Let $X \subseteq H^{*} = \sum_{i=1}^{n}\mathbb{Z}_{+} \epsilon_{i}$ and $X_{N} = X/NX$.

Define $\phi_{0} = \{\delta \in \Phi^{+}_{0} | \ \delta/2 \notin \Phi^{+}_{0} \}$ and
$\phi_{1} = \Phi_{1}^{+}$.
For a given $U_{q}^{(N)}(osp(1|2n))$ where $n \geq 2$ and $N \geq 3$ we define 
$\overline{\Lambda}_{N}^{+}$ by
$$\overline{\Lambda}_{N}^{+} = \left\{ \lambda \in X | \ 0 \leq \frac{2(\lambda + \rho,\alpha)}{(\alpha,\alpha)} \leq N', 
\ \forall \alpha \in \phi \right\},$$
where $\phi = \phi_{0} \cup \phi_{1}$ if $N \equiv 2 \pmod{4}$ and $\phi = \phi_{0}$ otherwise.
Here $N' = N$ if $N$ is odd and $N' = N/2$ otherwise.  
For $U_{q}^{(N)}(osp(1|2))$ where $N \geq 3$ we define
$\overline{\Lambda}_{N}^{+}$ by
$$\overline{\Lambda}_{N}^{+} = \left\{ \lambda \in X | \ 0 \leq (\lambda + \rho,\alpha) \leq N'' \right\},$$ where
$\alpha$ is the single odd root and $N'' = N/4$ if $N \equiv 2 \pmod{4}$ and $N'' = N'$ otherwise.
For each $\overline{\Lambda}_{N}^{+}$ we define $\Lambda_{N}^{+}$ identically except that we replace $\leq$ by $<$.  
$\Lambda_{N}^{+}$ plays the role of
the index set $I$ of the collection of $U_{q}^{(N)}(osp(1|2n))$ modules.

Let $V$ be the fundamental module of $U_{q}^{(N)}(osp(1|2n))$ with highest weight $\epsilon_{1}$. 
$V$ is irreducible, $2n+1$ dimensional and has the same structure as the fundamental module of $U_{q}(osp(1|2n))$ 
\cite{b}.  In \cite{b} we prove the following lemmas and theorems.
\begin{lemma}  
Set $N \geq 4$ to be even. 
For each $\mu \in \overline{\Lambda}_{N}^{+}$ there exists a finite dimensional 
$U_{q}^{(N)}(osp(1|2n))$ module $V_{\mu}$ with
highest weight $\mu$ such that $V_{\mu} \subseteq V^{\otimes t}$ for some $t \geq 1$.  
The quantum superdimension of $V_{\mu}$ is $sdim_{q}(V_{\mu}) = \lim_{q^{N} \rightarrow 1}(sdim_{q}(V_{\mu}^{gen}))$ \
where $V_{\mu}^{gen}$
is the finite dimensional irreducible $U_{q}(osp(1|2n))$ module with highest weight $\mu$.
The $sdim_{q}(V_{\mu}) \neq 0$ for all $\mu \in \Lambda_{N}^{+}$ and $sdim_{q}(V_{\mu}) = 0$ for all 
$\mu \in \overline{\Lambda}_{N}^{+} \backslash \Lambda_{N}^{+}$.
\end{lemma}

Let $W$ be the Weyl group of $osp(1|2n)$ and $\tau$ be the maximal element 
of $W$.  Then $-\tau(\lambda) \in \Lambda_{N}^{+}, \ \forall \lambda \in \Lambda_{N}^{+}$,
which implies that for each $V_{\lambda}$ there is a dual
module $(V_{\lambda})^{\dagger}$ with highest weight in $\Lambda_{N}^{+}$.

In theorems 2.1 and 2.2 we assume that $n \geq 2$ or that $n=1$ and $N \neq 4,8$.
\begin{theorem}
Set $N \geq 4$ to be even, $V$ to be the fundamental module of $U_{q}^{(N)}(osp(1|2n))$ 
and $\epsilon_{1} \in \Lambda_{N}^{+}$.  Then for every $t \in \mathbb{Z}_{+}$,
$$V^{\otimes t} = {\cal V} \oplus {\cal Z}$$ where
${\cal V} = \bigoplus_{\lambda \in \Lambda_{N}^{+}} \left(V_{\lambda} \right)^{\oplus m(\lambda)}$,
$m(\lambda) \in \mathbb{Z}_{+}$ and
${\cal Z}$ is a possibly empty direct sum of indecomposable modules, each with zero quantum superdimension.
\end{theorem}

\begin{theorem}
Set $N \geq 4$ to be even and let $V_{\lambda_{i}}$ be finite dimensional $U_{q}^{(N)}(osp(1|2n))$ modules with 
$\lambda_{i} \in \Lambda_{N}^{+}$ for all $i$.  Then
$$V_{\lambda_{1}} \otimes V_{\lambda_{2}} \otimes \cdots \otimes V_{\lambda_{s}} = {\cal V}' \oplus {\cal Z}'$$
where $s \geq 1$ and ${\cal V}'$ and ${\cal Z}'$ have the same form as in theorem 2.1
\end{theorem}

\begin{lemma}
Set $N \geq 3$ to be odd, $\epsilon_{1} \in \Lambda_{N}^{+}$ and $t \in \mathbb{Z}_{N/2+1/2-n}$.  Then
$V^{\otimes t} = \overline{{\cal V}}$ where $\overline{{\cal V}}$ has the same form as ${\cal V}$ does in theorem 2.1
\end{lemma}

\section{Finding the set $\{d_{\lambda}\}$}
We know that properties i), ii) and iii) hold for $U_{q}^{(N)}(osp(1|2n))$ at even $N$.  
We now consider condition
iv): there exists at least one collection of $\{d_{\lambda}\}$ that solves
$\chi_{\mu}(v) = \sum_{\lambda \in \Lambda_{N}^{+}}d_{\lambda}
\chi_{\lambda}(v^{-1})\chi_{\mu}(C_{\lambda})$
for all $\mu \in \Lambda_{N}^{+}$.  Now at generic $q$, 
the eigenvalue of $C_{\lambda}$ in an irreducible representation with highest weight $\mu \in \Lambda_{N}^{+}$ is given by
$sch_{\lambda}(q^{2(\mu + \rho)}) = S_{\lambda,\mu}/Q_{\mu}$ where 
$S_{\lambda,\mu} = (-1)^{[\lambda]}\sum_{\sigma \in W}\epsilon'(\sigma) q^{2(\lambda + \rho, \sigma(\mu + \rho))}$,
 $Q_{\mu} = \sum_{\sigma \in W}\epsilon'(\sigma)q^{2(\rho,\sigma(\mu + \rho))}$,
and where $\epsilon'(\sigma) = -1$ if the number of components of $\sigma$ 
that are reflections with respect to the elements of
$\phi_{0}$ is odd and $\epsilon'(\sigma) = 1$ otherwise.
Our proofs for theorems 2.1 and 2.2 tell us that $sch_{\lambda}(q^{2(\mu + \rho)})$ is well behaved when $q$ is taken to
the $N^{th}$ root of unity, and yields the desired $\chi_{\mu}(C_{\lambda})$.
To simplify finding the $\{d_{\lambda}\}$ we initially consider
$Q_{\mu} q^{-(\mu + 2\rho,\mu)} = 
\sum_{\lambda \in \Lambda_{N}^{+}} d_{\lambda}' q^{(\lambda +
2\rho,\lambda)}S_{\lambda,\mu}'$
where $d_{\lambda}' = (-1)^{[\lambda]}d_{\lambda}$
and $S_{\lambda,\mu}' = (-1)^{[\lambda]}S_{\lambda,\mu}$.
To solve for the $\{d_{\lambda}'\}$ we consider
$$Q_{\mu} q^{-(\mu + 2\rho,\mu)} = 
\sum_{\lambda \in X_{N}} x_{\lambda} q^{(\lambda +
2\rho,\lambda)}S_{\lambda,\mu}',$$
and set $x_{\lambda} = c q^{-(\lambda,2\rho)}$. We then have
$$Q_{\mu} q^{-(\mu + 2\rho,\mu)} = 
\sum_{\sigma \in W}\epsilon'(\sigma)
\sum_{\lambda \in X_{N}} cq^{(\lambda,\lambda)}
q^{2(\lambda + \rho,\sigma(\mu + \rho))}.$$
To ensure $x_{\lambda}$ is independent of $\mu$ we undertake the mapping 
$\lambda \rightarrow \sigma(\lambda + \rho) - \sigma(\mu + \rho)$ in the
summation, which may be done as $\sigma(\lambda + \rho) - \sigma(\mu + \rho) \in
X_{N}$ and the summation remains over $X_{N}$.
We then obtain
$$Q_{\mu}q^{-(\mu,\mu + 2\rho)} = \sum_{\lambda \in X_{N}}
cq^{(\lambda,\lambda + 2\rho) - (\mu,\mu + 2\rho)}\sum_{\sigma \in
W}\epsilon'(\sigma)q^{2(\rho,\sigma(\mu + \rho))}$$  which results in
$c^{-1} = \sum_{\lambda \in X_{N}} q^{(\lambda, \lambda + 2\rho)}$
and $$x_{\lambda} = q^{-(\lambda,2\rho)}/\sum_{\lambda \in
X_{N}}q^{(\lambda,\lambda + 2\rho)}.$$

Now $q^{(\lambda',\lambda' + 2\rho)} = (-1)^{p}q^{(\lambda,\lambda + 2\rho)}$
where $\lambda' = \lambda + N/2 \epsilon_{i}$ for any $\epsilon_{i}$ and $p = 1$ if $N = 4k$, $p = 0$ if $N =
2(2k+1)$.  Then $c$ is not well defined if $N = 4k$ and we set $N = 2(2k+1)$ 
for the remainder of this paper.  Now
$x_{\lambda'}q^{(\lambda',\lambda' + 2\rho)}S_{\lambda',\mu}' = x_{\lambda}q^{(\lambda,\lambda + 2\rho)}S_{\lambda,\mu}'$
where $\lambda' = \lambda + N/2 \epsilon_{i}$ for any $\epsilon_{i}$.  It then follows that
$\sum_{\lambda \in X_{N}} x_{\lambda} q^{(\lambda + 2\rho,\lambda)}S_{\lambda,\mu}'
= 2^{n}\sum_{\lambda \in X_{N/2}} x_{\lambda} q^{(\lambda + 2\rho,\lambda)}S_{\lambda,\mu}'.$

Let ${\overline{N}} = N/2$ and $\overline{\Lambda}_{\overline{N}}^{+}$ be the fundamental domain for 
$X_{\overline{N}}$ under
the action of the affine Weyl group $W_{\overline{N}}$ of $U_{q}^{(N)}(so_{2n+1})$ \cite{z1}.  
The affine Weyl groups $W_{N/2}$ of
$U_{q}^{(N)}(osp(1|2n))$ and $U_{q}^{(N)}(so_{2n+1})$ are identical and as 
$\overline{\Lambda}_{\overline{N}}^{+} = \overline{\Lambda}_{N}^{+}$
  it follows that $\overline{\Lambda}_{N}^{+}$ 
is a fundamental domain for $X_{N/2}$ under the action of $W_{N/2}$.

Now $S_{\sigma(\lambda + \rho) - \rho,\mu}'= \epsilon'(\sigma)S_{\lambda,\mu}'$ for any $\sigma \in W$.
If $\lambda \in \overline{\Lambda}_{N}^{+} \backslash \Lambda_{N}^{+}$  there either exists some $\sigma \in W$ such that
$\sigma(\lambda + \rho) - \rho = \lambda$ and $\epsilon'(\sigma) = -1$ or some $w \in W$ such that
$\epsilon'(w) = 1$ and $\lambda = w(\lambda + \rho) - \rho + k N/2 \epsilon_{i}$ for some $\epsilon_{i}$ and $k \in
\mathbb{Z}$.  
As $S_{\lambda',\mu}' = -S_{\lambda,\mu}'$ where $\lambda' = \lambda + N/2 \epsilon_{i}$, it follows
that $S_{\lambda,\mu}' = 0$ for $\lambda \in \overline{\Lambda}_{N}^{+} \backslash \Lambda_{N}^{+}$.

Then
$$2^{n} \sum_{\nu \in X_{N/2}} x_{\nu} q^{(\nu, \nu + 2\rho)}S_{\nu,\mu}'=
2^{n}\sum_{\lambda \in \Lambda_{N}^{+}} \sum_{\sigma \in W}\epsilon'(\sigma)x_{\sigma(\lambda + \rho)-\rho}
q^{(\lambda,\lambda + 2\rho)}S_{\lambda,\mu}'.$$
As $$\sum_{\lambda \in \Lambda_{N}^{+}} d_{\lambda}' q^{(\lambda +
2\rho,\lambda)}S_{\lambda,\mu}' = 2^{n}\sum_{\lambda \in \Lambda_{N}^{+}} \sum_{\sigma \in W}\epsilon'(\sigma)x_{\sigma(\lambda + \rho)-\rho}
q^{(\lambda,\lambda + 2\rho)}S_{\lambda,\mu}'$$ we obtain
$$d_{\lambda} = q^{(2\rho,\rho)}\gamma. sdim_{q}(V_{\lambda})/\sum_{\mu \in X_{N/2}}q^{(\mu, \mu + 2\rho)}$$
where $\gamma$ is the denominator in the expression of the quantum superdimension.  
Note that $d_{\lambda^{*}} = d_{\lambda}$ where
$\lambda^{*} = -\tau(\lambda)$.  The denominator of the $d_{\lambda}$ does not vanish: 
$\sum_{\mu \in X_{N/2}}q^{(\mu, \mu + 2\rho)}
= \sum_{\mu \in X_{N}}q^{(\mu, \mu + 2\rho)}/2^{n}$ and $\sum_{\mu \in X_{N}}q^{(\mu, \mu + 2\rho)}
= \prod_{k=0}^{n-1} G(N,2k+1)$ where $G(N,m) = \sum_{i=0}^{N-1}q^{i(i+m)}=
(1+i)\sqrt{N}/x^{m^{2}}$ and $x$ is a complex primitive $4N^{th}$ root of unity.

\end{section}

\begin{section}{Constructing the invariant}
The final matter we need to consider is condition v): that $$z = \sum_{\lambda \in \Lambda_{N}^{+}}d_{\lambda} q^{-(\lambda +
2\rho,\lambda)}sdim_{q}(V_{\lambda}) \neq 0.$$  It follows from \cite{tw} that $z \neq 0$.
Given $\zeta = \sum_{\lambda \in \Lambda_{N}^{+}} d_{\lambda} sdim_{q}(\lambda)$, $\zeta \neq 0$ and $z = d_{0} \zeta$.

Now it is a relatively simple matter to construct the invariants.  Denote a framed link in $S^{3}$ by $L$ and the
3-manifold it gives rise to by $M_{L}$.  Let $\sum(L)$ stand for the Reshetikhin-Turaev functor 
applied to $L$ (see \cite{rt,z1,z2}).  
Set $A_{L}$ to be the linking matrix of $L$ defined by: $a_{ii}$ is the framing number of the
$i^{th}$ component of $L$ and $a_{ij}$, $i \neq j$ is the linking number between the $i^{th}$ and $j^{th}$ components of 
$L$.  Let $\sigma(A_{L})$ be the number of nonpositive eigenvalues of $A_{L}$.  Then
$${\cal F}(M_{L}) = z^{-\sigma(A_{L})} \sum(L)$$
is a topological invariant of $M_{L}$.

\end{section}


\begin{thebibliography}{99}

\bibitem{rt}
Reshetikhin N Y and Turaev V G 1991
Invent. Math. 103 547-597

\bibitem{tw}
Turaev V G and Wenzl H 1993
Int. J. Math. 4
no. 2 323-358

\bibitem{z1} 
Zhang R B 1996
Comm. Math. Phys. 182 
no. 3 619-636

\bibitem{z2}
Zhang R B 
1994
Mod. Phys. Lett. A 9
no. 16 1453-1465,
1995
Rev. Math. Phys. 7
no. 5 809-831,
1997
Lett. Math. Phys. 41
no. 1 1-11

\bibitem{kt}
Khoroshkin S M and Tolstoy V N 1991
Comm. Math. Phys. 141
no. 3 599-617

\bibitem{z5}
Zhang R B 1992
J. Math. Phys. 33
no. 11 3918-3930

\bibitem{b}
Blumen S C
``Tensor products of finite dimensional $U_{q}(osp(1|2n))$ modules at roots of unity'', to appear.

\end{thebibliography}
\end{document}